\documentclass[11pt]{amsart}
\usepackage{CJK}
\usepackage{amssymb}
\usepackage{amsmath}
\usepackage{amsthm}

\usepackage{graphicx}
\usepackage{pst-plot}

\newtheorem{Prop}{Proposition}[section]
\newtheorem{Thm}[Prop]{Theorem}
\newtheorem{Lem}[Prop]{Lemma}
\newtheorem{Cor}[Prop]{Corollary}

\newtheorem{Q}[Prop]{Question}
\newtheorem{Fact}[Prop]{Fact}

\theoremstyle{definition}
\newtheorem{Ex}[Prop]{Example}

\begin{document}

\title{Semigroups of $L$-space Knots and Nonalgebraic Iterated Torus Knots}

\author{Shida Wang}

\address{Department of Mathematics, Indiana University, Bloomington, IN 47405}
\email{wang217@indiana.edu}

\begin{abstract}

Algebraic knots are known to be iterated torus knots and to admit\linebreak $L$-space surgeries.
However, Hedden proved that there are iterated torus knots that admit $L$-space surgeries but are not algebraic.
We present an infinite family of such examples, with the additional property that no nontrivial linear combination of knots in this family is concordant to a linear combination of algebraic knots.
The proof uses the Ozsv\'{a}th-Stipsicz-Szab\'{o} Upsilon function, and also introduces a new invariant of $L$-space knots, the formal semigroup.

\end{abstract}

\maketitle

\section{Introduction}

An algebraic knot $K$ can be defined to be the connected link of an isolated singularity of a complex curve in $\mathbb{C}^2$~\cite{Puiseux1,Puiseux2}.
All such knots are iterated torus knots but not vice versa~\cite[p.52]{Puiseux2} (we only consider positively iterated torus knots in this paper).
To be precise, an iterated torus knot $(((T_{p_1,q_1})_{p_2,q_2})\cdots)_{p_m,q_m}$ is an algebraic knot if and only if the indices satisfy $q_{i+1}>p_iq_ip_{i+1}$~\cite[Section 17a)]{Puiseux2}.

To each algebraic knot, one can associate a numerical semigroup of the nonnegative integers, denoted by $S_K$.
Initially this was done using the analytic properties of the curve, but $S_K$ is determined by the Alexander polynomial of $K$.
For example, for the torus knot $T_{p,q}$, $S_{T_{p,q}}=\langle p,q\rangle\subset\mathbb{Z}_{\geqslant0}$.
For algebraic knots, $S_K$ completely determines the Heegaard Floer complex $CFK^\infty(K)$.

By~\cite[Theorem 1.10]{cableL}, algebraic knots are all $L$-space knots, a class of knots defined using Heegaard Floer theory~\cite{Lknot}.
In this paper we will associate to each $L$-space knot what we call a \emph{formal} semigroup $S_K$, a subset of $\mathbb{Z}_{\geqslant0}$,
but now $S_K$ is not necessarily a semigroup.
Again, $S_K$ is determined by the Alexander polynomial of $K$ and it determines $CFK^\infty(K)$.
The fact that $S_K$ is not a semigroup provides an easily computed obstruction to $L$-space knots being algebraic.

We will show that for any $L$-space knot that is an iterated torus knot, $S_K$ is a semigroup.
However, we will use the Upsilon invariant recently defined by Ozsv\'{a}th, Stipsicz and Szab\'{o} in~\cite{upsilon} to show that many such $L$-space knots are not algebraic.
Going beyond this, we provide an infinite family of $L$-space iterated torus knots with the property that no nontrivial linear combination of these knots is even concordant to a connected sum of algebraic knots.
In particular, letting $\mathcal{C}$ denote the smooth concordance group and $\mathcal{C}_A$ the subgroup generated by algebraic knots, we prove the following:

\begin{Thm}\label{main}$\mathcal{C}/\mathcal{C}_A$ is infinitely generated.\end{Thm}

Note that $\mathcal{C}_A$ is also infinitely generated, even restricted to algebraically slice\linebreak knots~\cite{infinite-algebraic}.

We will compute the $\Upsilon$ functions of this infinite family of $L$-space iterated torus knots and prove they cannot be generated by $\Upsilon$ functions of $(n,n+1)$-torus knots.
Hence the following result of Feller and Krcatovitch implies Theorem~\ref{main}.

\begin{Thm}\label{expansion}\emph{(\cite[Proposition 6 and the paragraph before it]{key})} The $\Upsilon$ function of any algebraic knot is a sum of $\Upsilon$ functions of $(n,n+1)$-torus knots.\end{Thm}

In the computation, we will observe the behavior of $S_K$ for $L$-space knots under cabling operation (see Proposition~\ref{algorithm}).

Hedden proved that if $K$ is an $L$-space knot and $q\geqslant p(2g(K)-1)$, then the cable $K_{p,q}$ is an $L$-space knot~\cite[Theorem 1.10]{cableL}. In~\cite{cableLconverse}, Hom proved that the converse is true.

\begin{Thm}\label{cableLiff}Assume that $K\subset S^3$ is a nontrivial knot and $p\geqslant2$. The $(p,q)$-cable of a knot $K$ is an $L$-space knot if and only if $K$ is an $L$-space knot and $q\geqslant p(2g(K)-1)$.\end{Thm}

We will prove an analogue of this statement.

\begin{Thm}Assume that $K\subset S^3$ is a nontrivial knot and $p\geqslant2$.
The $(p,q)$-cable of a knot $K$ is an $L$-space knot with $S_{K_{p,q}}$ being a semigroup if and only if $K$ is an $L$-space knot with $S_K$ being a semigroup and $q\geqslant p(2g(K)-1)$.\end{Thm}

\emph{Acknowledgments.} The author wishes to express sincere thanks to Professor Charles Livingston for proposing this study and carefully reading a draft of this paper.
Thanks also to David Krcatovich for a motivating conversation in December 2015 and sharing his computer programs
and to Professor Yi Ni for an email correspondence about Question~\ref{semigroupConj}.

\section{Formal Semigroups Under Cabling}\label{SectionSemigroup}

\subsection{Formal semigroups of $L$-space knots}

Write $\mathbb{Z}_{>k}:=\{m\in\mathbb{Z}|m>k\}$ and
$\mathbb{Z}_{\geqslant k}:=\{m\in\mathbb{Z}|m\geqslant k\}$.

For any $L$-space knot $K$, the Alexander polynomial $\Delta_K(t)=\sum_{i=0}^{2n}(-1)^it^{\alpha_i}$,
where $0=\alpha_0<\alpha_1<\cdots<\alpha_{2n}$~\cite[Theorem 1.2]{Lknot} and $\frac{\alpha_{2n}}{2}=g(K)$ is the genus of $K$~\cite[Theorem 1.2]{genus}.
We will not use the symmetrized Alexander polynomial.

Consider $\Delta_K(t)$ as an element in the ring $\mathbb{Z}[[t]]$ of formal power series with integer coefficients.
Define the \emph{formal semigroup} $S_K$ of the $L$-space knot $K$ to be the subset of $\mathbb{Z}_{\geqslant0}$ satisfying $\sum_{s\in S_K}t^s=\frac{\Delta_K(t)}{1-t}$,
where the right-hand side is sometimes called the Alexander function.
Since $\Delta_K(t)=\sum_{i=0}^{2n}(-1)^it^{\alpha_i}$, it follows that $$S_K=\{\alpha_0,\cdots,\alpha_1-1,\alpha_2,\cdots,\alpha_3-1,\cdots,\alpha_{2n-2},\cdots,\alpha_{2n-1}-1,\alpha_{2n}\}\cup\mathbb{Z}_{>\alpha_{2n}}.$$

\emph{Remark.} (i) $S_K$ is denoted $\Gamma_K$ in \cite{enumerating}.

(ii) $\alpha_1=1$. More generally, the $(i+1)$th element in $S_K$ is bounded below by $2i$ for $0\leqslant i\leqslant g(K)$~\cite[Theorem 1.6]{geography}.

\begin{Ex}Let $K$ be the torus knot $T_{3,7}$.

$$\Delta_K(t)=\frac{(t^{21}-1)(t-1)}{(t^3-1)(t^7-1)}=1-t+t^3-t^4+t^6-t^8+t^9-t^{11}+t^{12}$$ $$=(1-t)(1+t^3+t^6+t^7+t^9+t^{10}+t^{12}+\sum_{s>12}t^s).$$

So $S_K=\{0,3,6,7,9,10,12\}\cup\mathbb{Z}_{>12}=\langle 3,7\rangle$.
\end{Ex}

\begin{Fact}\emph{(\cite{Puiseux1})} For algebraic knots, $S_K$ is a semigroup, and it equals the traditionally defined semigroup of the link of singularity.\end{Fact}

\emph{Remark.} No matter whether $S_K$ is a semigroup, it is dual with respect to $2g(K)-1$.
That is, $s\in S_K\Leftrightarrow2g(K)-1-s\not\in S_K$.
This follows from the palindromicity of the symmetrized Alexander polynomial.

\begin{Ex}\label{notSemigroup}(\cite[Example 2.3]{enumerating}) The pretzel knot $P(-2,3,7)$ has $S_K=\{0,3,5,7,8,10\}\cup\mathbb{Z}_{>10}$, which is not a semigroup.

Generally, for any odd integer $n\geqslant7$ the pretzel knot $P(-2,3,n)$ is an $L$-space knot~\cite{Lknot}.
By a recursive formula for the Alexander polynomial of $(-2,3,n)$ pretzel knots (cf.~\cite[Equation (1---3)]{pretzel}), one can verify $S_{P(-2,3,n)}\cap[0,7]=\{0,3,5,7\}$ for any $n$.
So $S_{P(-2,3,n)}$ is not a semigroup.\end{Ex}

\subsection{The cabling formula}

Let $K$ be a nontrivial $L$-space knot.
We will give a formula in Proposition~\ref{algorithm} and use it to prove the following statement.

\begin{Thm}\label{cableS}Let $K$ be a nontrivial $L$-space knot and $q\geqslant p(2g(K)-1)$.
Then $S_K$ is a semigroup if and only if $S_{K_{p,q}}$ is a semigroup.\end{Thm}

Then it is easy to show that

\begin{Cor}\label{iteratedS}If an $L$-space knot $K$ is an iterated torus knot, then $S_K$ is a semigroup.\end{Cor}

\begin{Ex}Let $K=(T_{2,3})_{2,k}$ where $k$ is an odd integer.
Then $K$ is an $L$-space knot if $k\geqslant3$.
Also, $K$ is an algebraic knot if and only if $k\geqslant13$~\cite[Section 17a)]{Puiseux2}.
So if $3\leqslant k<13$, then $K$ is not an algebraic knot but $S_K$ is still a semigroup.
\end{Ex}

Theorem \ref{cableS} is based on the following fact.

\begin{Prop}[Cabling formula]\label{algorithm}Let $K$ be a nontrivial $L$-space knot. Suppose $p\geqslant2$ and $q\geqslant p(2g(K)-1)$.
Then $S_{K_{p,q}}=pS_K+q\mathbb{Z}_{\geqslant0}:=\{pa+qb|a\in S_K,b\in\mathbb{Z}_{\geqslant0}\}$.\end{Prop}

\textbf{Proof.} Recall that $\Delta_{K_{p,q}}(t)=\Delta_K(t^p)\Delta_{T_{p,q}}(t)$.

So $\frac{\textstyle\Delta_{K_{p,q}}(t)}{\textstyle1-t}=\frac{\textstyle\Delta_K(t^p)}{\textstyle1-t}\cdot\frac{\textstyle(t^{pq}-1)(t-1)}{\textstyle(t^p-1)(t^q-1)}
=\frac{\textstyle\Delta_K(t^p)}{\textstyle1-t^p}\cdot\frac{\textstyle t^{pq}-1}{\textstyle t^q-1}$.

By definition $\sum_{s\in S_K}t^s=\frac{\textstyle \Delta_K(t)}{\textstyle1-t}$.
Hence $\sum_{s\in pS_K}t^s=\frac{\textstyle\Delta_K(t^p)}{\textstyle1-t^p}$.
Observe that $\frac{\textstyle t^{pq}-1}{\textstyle t^q-1}=1+t^q+\cdots+t^{(p-1)q}$.
Therefore$$\frac{\textstyle\Delta_{K_{p,q}}(t)}{\textstyle1-t}=(\sum_{s\in pS_K}t^s)\cdot(1+t^q+\cdots+t^{(p-1)q}).$$
By definition $\sum_{s\in S_{K_{p,q}}}t^s=(\sum_{s\in pS_K}t^s)\cdot(1+t^q+\cdots+t^{(p-1)q})$.

Now $(\sum_{s\in pS_K}t^s)\cdot(1+t^q+\cdots+t^{(p-1)q})\\=\sum_{s\in pS_K}t^s+\sum_{s\in pS_K}t^{s+q}+\cdots+\sum_{s\in pS_K}t^{s+(p-1)q}$.

To show $S_{K_{p,q}}=pS_K+q\mathbb{Z}_{\geqslant0}$, it suffices to prove that $pS_K+q\mathbb{Z}_{\geqslant0}$ is the disjoint union of $pS_K,\;pS_K+q,\;\cdots,\;pS_K+(p-1)q$.

The sets $pS_K,\;pS_K+q,\;\cdots,\;pS_K+(p-1)q$ must be pairwise disjoint. Otherwise some term of $\sum_{s\in S_{K_{p,q}}}t^s$ would have coefficient greater than 1.

Next, $(pS_K)\cup(pS_K+q)\cup\cdots\cup(pS_K+(p-1)q)\subset pS_K+q\mathbb{Z}_{\geqslant0}$ clearly.

To prove $(pS_K)\cup(pS_K+q)\cup\cdots\cup(pS_K+(p-1)q)\supset pS_K+q\mathbb{Z}_{\geqslant0}$, let\linebreak $pa+qb\in S_K+q\mathbb{Z}_{\geqslant0}$, where $a\in S_K,b\in\mathbb{Z}_{\geqslant0}$.
Suppose $b=kp+c$ with $k\in\mathbb{Z}_{\geqslant0}$ and $c\in\{0,1,\cdots,p-1\}$.
Then $pa+qb=pa+q(kp+c)=p(a+kq)+cq$. It suffices to show $p(a+kq)\in pS_K$.
If $k=0$, this is trivial. If $k>0$, then $a+kq\geqslant q\geqslant p(2g(K)-1)\geqslant2g(K)$, since we assumed $p\geqslant2$.
Hence $a+kq\in S_K$ by the fact that $\mathbb{Z}_{\geqslant2g(K)}\subset S_K$.
\hfill$\Box$

\textbf{Proof of Theorem~\ref{cableS}.} The proof in the case of $p=1$ is trivial. Assume $p\geqslant2$.

If $S_K$ is a semigroup, then $S_{K_{p,q}}=pS_K+q\mathbb{Z}_{\geqslant0}$ is a semigroup.

If $S_K$ is not a semigroup, then since $\mathbb{Z}_{\geqslant2g(K)}\subset S_K$, there are $x,y\in S_K$ such that $x+y\not\in S_K$ and $x+y<2g(K)$.
So $px,py\in S_{K_{p,q}}$. It suffices to show $px+py\not\in S_{K_{p,q}}$.
Observe that $px+py=p(x+y)\leqslant p(2g(K)-1)<q$, where $p(2g(K)-1)\neq q$ because $p$ and $q$ are relatively prime.
Thus, if $px+py=pa+qb$ for some $a\in S_K,b\in\mathbb{Z}_{\geqslant0}$, then $b$ must be $0$.
Therefore $px+py=pa\Rightarrow x+y=a\in S_K$, which is impossible.
\hfill$\Box$

\textbf{Proof of Corollary~\ref{iteratedS}.} Suppose $(((T_{p_1,q_1})_{p_2,q_2})\cdots)_{p_m,q_m}$ is an $L$-space knot.
Then $(((T_{p_1,q_1})_{p_2,q_2})\cdots)_{p_k,q_k}$ is an $L$-space knot for $k=2,\cdots,m$ and\\
$q_k\geqslant p_k(2g((((T_{p_1,q_1})_{p_2,q_2})\cdots)_{p_{k-1},q_{k-1}})-1)$ by Theorem~\ref{cableLiff}.
Hence the conclusion follows from Theorem~\ref{cableS}.
\hfill$\Box$

\emph{Remark.} In fact, Proposition~\ref{algorithm} gives an algorithm to compute generators of $S_K$ for $K=(((T_{p_1,q_1})_{p_2,q_2})\cdots)_{p_m,q_m}$.
A set of generators is $$\{p_1p_2\cdots p_m,\ q_1p_2\cdots p_m,\ q_2p_3\cdots p_m,\ \cdots,\ q_{m-1}p_m,\ q_m\}.$$

It is natural to ask the following question.

\begin{Q}\label{semigroupConj}Is there an $L$-space knot $K$ with $S_K$ being a semigroup, but $K$ is not an iterated torus knot?\end{Q}

Similarly to the motivation of~\cite[Conjecture 1.3]{surgery}, if the answer is ``no'', then the surgery coefficient of any finite surgery on any hyperbolic knot must be an integer by~\cite[Theorem 1.2]{surgery}.

The author did not find any examples for a ``yes'' answer by computing Alexander polynomials for some $L$-space knots provided in~\cite{TT} and~\cite{Lexamples}.

\section{A Family of Nonalgebraic $L$-space Iterated Torus knots}

The papers \cite{deformation,enumerating,rational} do not use the fact that $S_K$ is a semigroup, so their conclusions are true for $L$-space knots and Section~\ref{SectionSemigroup} does not apply.
However, \cite{cobordism} uses the fact that $S_K$ is a semigroup, so its conclusion can be generalized.

\subsection{Review of the Upsilon invariant}

We refer to \cite{upsilon} for the definition of the Upsilon invariant.
For our purpose, we only need to know the following properties.

\begin{Thm}\label{properties}\emph{(\cite[Section 1]{upsilon})} For each $t\in[0,2]$ there is a well-defined knot invariant $\Upsilon_K(t)$. Moreover, $\Upsilon_K(t)$ satisfies the following properties:

(i) $\Upsilon_K(t)$ is a piecewise linear function in $t$ on $[0,2]$.

(ii) $\Upsilon_K(t)=\Upsilon_K(2-t)$.

(iii) $\Upsilon_{-K}(t)=-\Upsilon_K(t)$ and $\Upsilon_{K_1\#K_2}(t)=\Upsilon_{K_1}(t)+\Upsilon_{K_2}(t)$.

(iv) $\Upsilon_K(t)=0$ if $K$ is smoothly slice.

(v) $\frac{t_0}{2}\Delta\Upsilon'_K(t_0)$ is an integer for any $t_0\in(0,2)$, where\\\hspace*{2em}$\Upsilon'_K(t_0):=\lim\limits_{t\rightarrow t_0+}\Upsilon_K'(t)-\lim\limits_{t\rightarrow t_0-}\Upsilon_K'(t)$.\end{Thm}

In \cite[Theorem 6.2]{upsilon}, the Upsilon invariant of $L$-space knots is computed in terms of the Alexander polynomial.

Alternatively, the Upsilon invariant can be expressed in terms of formal semigroups for $L$-space knots as follows, which was first stated in \cite[Proposition 4.4]{deformation} for algebraic knots.

\begin{Prop}\label{formula} Let $K$ be an $L$-space knot with genus $g$ and $S$ be the corresponding formal semigroup.
Then for any $t\in[0,2]$ we have $$\Upsilon_K(t)=\max\limits_{m\in\{0,\cdots,2g\}}\{-2\#(S\cap[0,m))-t(g-m)\}.$$\end{Prop}

The location of the first singularity (the discontinuity of the derivative) of the\linebreak Upsilon invariant for algebraic knots is given in \cite[Theorem 8]{cobordism}. This can be easily generalized to $L$-space knots with semigroups.

\begin{Thm}\label{1singularity}Let $K$ be an $L$-space knot with genus $g$.
If $S_K$ is a semigroup and the least nonzero element of $S_K$ is $a$, then $\Upsilon_K(t)=-gt$ for $t\in[0,\frac{2}{a}]$ and $\Upsilon_K(t)>-gt$ for $t>\frac{2}{a}$.\end{Thm}

\subsection{Upsilon invariant of algebraic knots}

\begin{Prop}Let $f(t)$ be a linear combination $\sum c_i\Upsilon_{T_{n_i,n_i+1}}(t)$ where $c_i\in\mathbb{Z}$.
Then $\Delta f'(\frac{2}{p})=\Delta f'(\frac{4}{p})$ for any odd integer $p\geqslant3$.\end{Prop}

\textbf{Proof.} Let $n$ be any positive integer.
According to \cite[Proposition 6.3]{upsilon},
$$\Delta\Upsilon'_{T_{n,n+1}}(t)=\left\{\begin{aligned}&n&\text{ for }t=\frac{2i}{n},0<i<n\\&0&\text{ otherwise}.\end{aligned}\right.$$
If $n$ does not divide $p$, then neither $\frac{2}{p}$ nor $\frac{4}{p}$ belongs to the set $\{\frac{2i}{n}\ \big|\ 0<i<n\}$.
Hence $\Delta\Upsilon'_{T_{n,n+1}}(\frac{2}{p})=\Delta\Upsilon'_{T_{n,n+1}}(\frac{4}{p})=0$.
If $p=kn$ for some $k\in\mathbb{Z}_{>0}$, then both $\frac{2}{p}$ and $\frac{4}{p}$ belong to the set $\{\frac{2i}{n}\ \big|\ 0<i<n\}$.
Hence $\Delta\Upsilon'_{T_{n,n+1}}(\frac{2}{p})=\Delta\Upsilon'_{T_{n,n+1}}(\frac{4}{p})=n$.

The conclusion follows from the fact that $f(t)$ is a linear combination\linebreak $\sum c_i\Upsilon_{T_{n_i,n_i+1}}(t)$.
\hfill$\Box$

Using Theorem~\ref{expansion}, we immediately obtain the following corollary.

\begin{Cor}\label{1equals2}If $K$ is an algebraic knot, then $\Delta\Upsilon'_K(\frac{2}{p})=\Delta\Upsilon'_K(\frac{4}{p})$ for any odd integer $p\geqslant3$.\end{Cor}

\subsection{A family of nonalgebraic knots}

Now we will consider the family of knots $\{J_k\}_{k=3}^\infty$ where $J_k=(T_{2,3})_{k,2k-1}$.
By Proposition~\ref{algorithm}, the formal semigroup\linebreak $S_{J_k}=\langle2k-1,\ 2k,\ 3k\rangle$.
The following corollary is an easy consequence of\linebreak Theorem~\ref{1singularity}.

\begin{Cor}\label{examples1}$\Upsilon_{J_k}(t)=-g(J_k)\:t$ for $t\in[0,\frac{2}{2k-1}]$ and\\\hspace*{7em} $\Upsilon_{J_k}(t)>-g(J_k)\:t$ for $t>\frac{2}{2k-1}$.\end{Cor}

The first singularity of $\Upsilon_{J_k}(t)$ is at $t=\frac{2}{2k-1}$. We will show that the second singularity is at $t=\frac{4}{k+1}$.

\begin{Lem}\label{examples2}$\Upsilon_{J_k}(t)=-2-(g(J_k)-(2k-1))\:t$ for $t\in[\frac{2}{2k-1},\frac{4}{k+1}]$ and\\\hspace*{6em} $\Upsilon_{J_k}(t)\geqslant-6-(g(J_k)-3k)\:t$ for $t\geqslant\frac{4}{k+1}$.\end{Lem}

\textbf{Proof.} Fix the integer $k\geqslant3$. Abbreviate $g(J_k)=g,\ S_{J_k}=S,\ \Upsilon_{J_k}=\Upsilon$.

Taking $m=2k-1$, we have the linear function $$-2\#(S\cap[0,m))-t(g-m)=-2-(g-(2k-1))\,t.$$ So $\Upsilon(t)\geqslant-2-(g-(2k-1))\,t$.

To show $\Upsilon_{J_k}(t)\leqslant-2-(g-(2k-1))\,t$ on $[\frac{2}{2k-1},\frac{4}{k+1}]$, we will consider the cases of $m=0$, \ $0<m\leqslant2k-1$, \ $m=2k$ and $m>2k$ separately.

If $m=0$, then $-2\#(S\cap[0,m))-t(g-m)$\\\hspace*{\fill}$=-g\,t\leqslant-g\,t+(2k-1)t-2=-2-(g-(2k-1))\,t$ since $t\leqslant\frac{2}{2k-1}$.

If $0<m\leqslant2k-1$, then $-2\#(S\cap[0,m))-t(g-m)$\\\hspace*{\fill}$=-2-(g-(2k-1))\,t$ since $t\leqslant\frac{2}{2k-1}$.

If $m=2k$, then $-2\#(S\cap[0,m))-t(g-m)=-4-(g-2k)\,t$\\\hspace*{\fill}$=-2-(g-(2k-1))\,t-2+2t\leqslant-2-(g-(2k-1))\,t$ since $t\leqslant2$.

If $m>2k$, this final case is the most delicate one. Here are the details.

We claim that $(k+1)\,(\#(S\cap[0,m))-1)\geqslant2(m-(2k-1))$.

This inequality can be simply verified as $(k+1)\,(3-1)\geqslant2(m-(2k-1))$ when $m\leqslant3k=2k-1+(k+1)$.
Without loss of generality, assume there is a positive integer $n$ such that $2k-1+n(k+1)<m\leqslant2k-1+(n+1)(k+1)$.
Since $S$ is generated by $2k-1$, $2k$ and $3k$, we have $0,\ 2k-1,\ 2k,\ 3k\in S$\linebreak and therefore $4k-1,\ 4k,\ 5k-1,\ 5k,\ 6k-1,\ 6k,\ \cdots\in S$.
Clearly\linebreak $0,2k-1,2k,3k,4k-2,4k-1,4k,5k-1,5k,\cdots,(2+n)k-1,(2+n)k\in S\cap[0,m)$.
Thus $\#(S\cap[0,m))\geqslant2(n+1)+1$ and therefore $$(k+1)\,(\#(S\cap[0,m))-1)\geqslant(k+1)(2(n+1)+1-1)$$ $$=2(2k-1+(n+1)(k+1)-(2k-1))\geqslant2(m-(2k-1)).$$

The claim implies $-2\#(S\cap[0,m))-t(g-m)\leqslant-2(\frac{2(m-(2k-1))}{k+1}+1)-tg+tm\\\leqslant\frac{-4(m-(2k-1))}{k+1}-2-gt+\frac{4}{k+1}m
=\frac{4(2k-1)}{k+1}-2-(g-(2k-1))\,t-(2k-1)t\\\leqslant\frac{4(2k-1)}{k+1}-2-(g-(2k-1))\,t-(2k-1)\frac{4}{k+1}=-2-(g-(2k-1))\,t$ since $t\leqslant\frac{4}{k+1}$.

To prove the second part of the lemma, take $m=3k$. Then $$-2\#(S\cap[0,m))-t(g-m)=-6-(g-3k)\,t.$$ So $\Upsilon(t)\geqslant-6-(g-3k)\,t$.
\hfill$\Box$

\begin{Thm}Let $\mathcal{C}_A$ be the subgroup of $\mathcal{C}$ generated by algebraic knots
and $\mathcal{G}$ be any subgroup of $\mathcal{C}$ such that $\mathcal{C}_A\subset\mathcal{G}$ and $J_k\in\mathcal{G},\forall k\geqslant3$.
Then $\{J_k\}_{k=3}^\infty$ generates a $\mathbb{Z}^\infty$ direct summand of $\mathcal{G}/\mathcal{C}_A$.\end{Thm}

\textbf{Proof.} By Theorem~\ref{properties}(v) and Corollary~\ref{1equals2}, we know that
$$\lambda_k:K\mapsto\frac{1}{2k-1}\Delta\Upsilon'_K(\frac{2}{2k-1})-\frac{1}{2k-1}\Delta\Upsilon'_K(\frac{4}{2k-1})$$
is a well-defined homomorphism from $\mathcal{G}/\mathcal{C}_A$ to $\mathbb{Z}$ for any integer $k\geqslant2$.
By Corollary~\ref{examples1} and Lemma~\ref{examples2}, we know that $\lambda_k(J_k)=1$ for any integer $k\geqslant3$.
Also, $\lambda_i(J_k)=0,\forall i>k$ by Lemma~\ref{examples2} and Proposition~\ref{1equals2}.
Hence $\{J_k\}_{k=3}^\infty$ generates a $\mathbb{Z}^\infty$ direct summand of $\mathcal{G}/\mathcal{C}_A$ by \cite[Lemma 6.4]{upsilon}.
\hfill$\Box$

Summarizing, we have:

\hspace*{2ex}\{algebraic knots\}\\
$\subset$\{$L$-space iterated torus knots\}\\
$\subset$\{$L$-space knots whose formal semigroup is a semigroup\}\hfil (by Corollary~\ref{iteratedS})\\
$\subset$\{$L$-space knots\}.

The knots $\{J_k\}$ lie in the first gap. Question~\ref{semigroupConj} asks whether the second gap is empty. Knots in Example~\ref{notSemigroup} lie in the third gap.

\end{document}